\documentclass[lettersize,journal]{IEEEtran}
\usepackage{amsmath,amsfonts}
\usepackage{subfigure}
\usepackage{algorithmic}
\usepackage{algorithm}
\usepackage{array}
\usepackage[caption=false,font=normalsize,labelfont=sf,textfont=sf]{subfig}
\usepackage{textcomp}
\usepackage{stfloats}
\usepackage{url}
\usepackage{verbatim}
\usepackage{graphicx}
\usepackage{cite}
\usepackage{amsmath}
\usepackage{subfig}
\usepackage{makecell}
\begin{document}
% ===============================================================================
%% Paper Title + etc
\title{Efficient Estimation of Active Element Patterns for 2-D Planar Array Antennas via Directional Decomposition}
% ===============================================================================
\author{Jeong-Wan Lee, Sung-Jun Yang,~\IEEEmembership{Member, IEEE}%
\thanks{This paper was financially supported by Seoul National University of Science and Technology. (\textit{Corresponding Author: Sung-Jun Yang.})}%
\thanks{The authors are with the Department of Electronic Engineering, Seoul National University of Science and Technology, Rebublic of South Korea (e-mail: leejeongwan@seoultech.ac.kr; wnstjddid@seoultech.ac.kr)}%

\thanks{Manuscript received February XX, 2025; revised February XX, 2025.}}%
\markboth{IEEE XX,~Vol.~XX, No.~X, February~2025}%
{Shell \MakeLowercase{\textit{et al.}}: A Sample Article Using IEEEtran.cls for IEEE Journals}
\maketitle
\begin{abstract}
The active element pattern method is widely employed in beam pattern synthesis of array antenna to account for mutual coupling between antenna elements.
Calculating the active element patterns for large number of array requires full-wave analyses of total array structure, which is time consuming.
To obtain accurate active element patterns efficiently, this letter proposes a method to estimates active element patterns in largely arrayed antenna using directional decomposition approach.
Reducing computational cost, proposed method constructs the transfer matrices to reflect both mutual coupling and truncation effects between each antenna element.
Numerical validation with open-ended waveguides confirms that the proposed method can estimate active element patterns with high accuracy.
The synthesized beam patterns show mean squared errors below 0.1dB in the main lobe region for various beam steering cases.
The computational complexity for numerical analysis reduces from $\mathcal{O}(M_B^2(N_x^3 N_y^3))$ to $\mathcal{O}(M_B^2(N_x^3 + N_y^3))$, resulting in a reduction of computation time to under 0.095\% compared to the conventional active element pattern method.

\end{abstract}

% ===============================================================================
%% 1. Introduction
﻿\section{Introduction}
\IEEEPARstart{A}{rray} antenna radiation pattern analysis plays a critical role in determining the electromagnetic performance of modern communication and radar systems. By arraying antennas, the electromagnetic interaction between antenna elements induces mutual coupling effects causing the radiation patterns to differ for each antenna element. The analysis of these effects requires significant computational costs, particularly for large two-dimensional (2-D) array antennas.

The pattern multiplication method (PMM) assumes that all array antenna elements share the same radiation pattern [4–8]. This method achieves minimal computational cost by requiring numerical analysis of only a single element’s radiation pattern, although the error increases as the mutual coupling effect becomes stronger [4–6]. For uniform array structures, a numerical technique applying periodic boundary conditions can partially reflect the mutual coupling effect between array elements [7–8]. However, this approach cannot account for the truncation effect at finite array edges due to its infinite array assumption.

To accurately model mutual coupling and truncation effects, full-wave analysis methods such as method of moment (MoM) and finite element method can be employed. Such approach is based on the active element pattern (AEP) method [1–3], which calculates radiation patterns by exciting only one antenna port while terminating the others. As the number of array elements increases, the costs to compute Apng for all ports increases rapidly, motivating various approaches to reduce computational demands [9–14]. Recently, machine learning–based methods demonstrated high accuracy in estimating far‐field Apng [12–14], but require additional full‐wave analyses for extensive training datasets. The coupling matrix approach has also been adopted for array radiation pattern analysis due to its intuitive characteristics [9-10]. This approach constructs mutual coupling matrices through open‐circuit analysis using the $Z$‐matrix of the array antenna [15–17]. However, most of these approaches still require full‐wave analysis for the entire structure of large 2‐D arrays.

The finite-by-infinite array analysis approach [18] was proposed to improve the efficiency by decomposing 2-D array problems into a series of one-dimensional (1-D) analyses. Using Floquet mode analysis, the approach demonstrates accurate analysis for mutual coupling effects. Despite this advantage, the infinite array assumption of Floquet mode limits its accuracy for edge elements, and the requirement for independent infinite‐array simulations for varying observation angles leads to additional computational loads for large 2‐D arrays with wide scan angles.
To overcome these limitations while preserving the computational advantages of directional decomposition approach, this letter presents a simple and efficient method for AEP estimation in large 2‐D planar arrays based on transfer matrix. Both mutual coupling and truncation effects on array structures can be accurately considered in the proposed approach. The proposed method employs current‐based transfer matrices which eliminates the need for additional Floquet mode analysis, replacing complex 2‐D array analysis with simpler 1‐D array problems. Numerical validation demonstrates that the proposed approach maintains high accuracy while decreasing computational time and complexity compared to conventional approaches.

This letter is organized as follows. Section II presents a brief background of the AEP method and describes a directional decomposition‐based AEP estimation method. Section III validates the performance of the proposed method through comparisons of results with conventional full‐wave AEP methods and PMM. Finally, Section IV concludes the letter.

% ===============================================================================
%% 2. Array Pattern Analysis with Active Element Pattern Method
\section{AEP Estimation via Directional Decomposition}
Accurate estimation of Apng requires consideration of both mutual coupling and truncation effects in array antenna structures.
This section presents a directional decomposition approach to AEP calculation using current-based transfer matrices for array beam synthesis.
Following the theoretical formulation of AEP method, complex 2-D array structures are decomposed into simplified 1-D array configurations.

\subsection{AEP Method for Array Beam Synthesis}
The AEP method characterizes individual element radiation in the presence of environmental effects of the array. For array beam synthesis, Apng provide accurate radiation patterns by accounting for both mutual coupling and edge truncation effects. 
Considering a 2-D array of ${N_x}\times{N_y}$ elements, radiated E-field array pattern of a finite array \(F^{AEP}(\theta, \phi)\) is calculated as
\begin{equation}
F^{AEP}(\theta, \phi) = \sum_{k=1}^{{N_x}{N_y}} w_k E^{AEP}_k(\theta, \phi)
\end{equation}
where \( w_k \) (\( k = 1, \ldots, {N_x} {N_y} \)) is arbitrary excitations of \(k\)-th port and \(E^{AEP}_k(\theta, \phi)\) indicates AEP (E-field) when \(k\)-th port is excited.
According to MoM formulation, $E^{AEP}_k(\theta, \phi)$ can be obtained as
\begin{equation}
   E^{\text{AEP}}_{k}(\theta,\phi) = \int_S   G(\mathbf{u}(\theta,\phi), \mathbf{r}') \cdot J^{\text{AEP}}_{k}(\mathbf{r}') \, d\mathbf{s}'
\end{equation}
where $\mathbf{u}(\theta,\phi)$ is the unit vector pointing in $\theta$ and $\phi$ direction, $\mathbf{r}'$ is current source point on antenna surface in area of surface $S$. $G(\mathbf{u}(\theta,\phi), \mathbf{r}')$ used in calculating Apng is free-space dyadic Green's function.
$J^{AEP}_{k}\in\mathbb{C}^{{N_x}{N_y}M}$ represents the induced current density distribution of array structure when \(k\)-th port is excited. The current distribution can be determined numerically by MoM for each antenna element segmented into M triangular meshes. As the number of arrays increases, it requires large computation time to obtain Apng due to the raised number of segmentation.
Therefore, the following subsection proposes an efficient method for estimating Apng through current-based transfer matrices, reducing computational cost while maintaining accuracy.

\subsection{Determination of Transfer Matrix on 1-D Array}
The proposed directional decomposition approach enables dimensional expansion of induced current distributions from 1-D to 2-D array structures. The Apng can then be calculated directly applying free-space dyadic Green's function to the current distributions of 2-D arrays. Unlike conventional coupling matrix methods [15-17], this approach determines the transfer matrix based on current distributions in 1-D array configurations, as illustrated in Fig. 1.

\begin{figure}[h]
 \centering
 \includegraphics[width=0.48\textwidth]{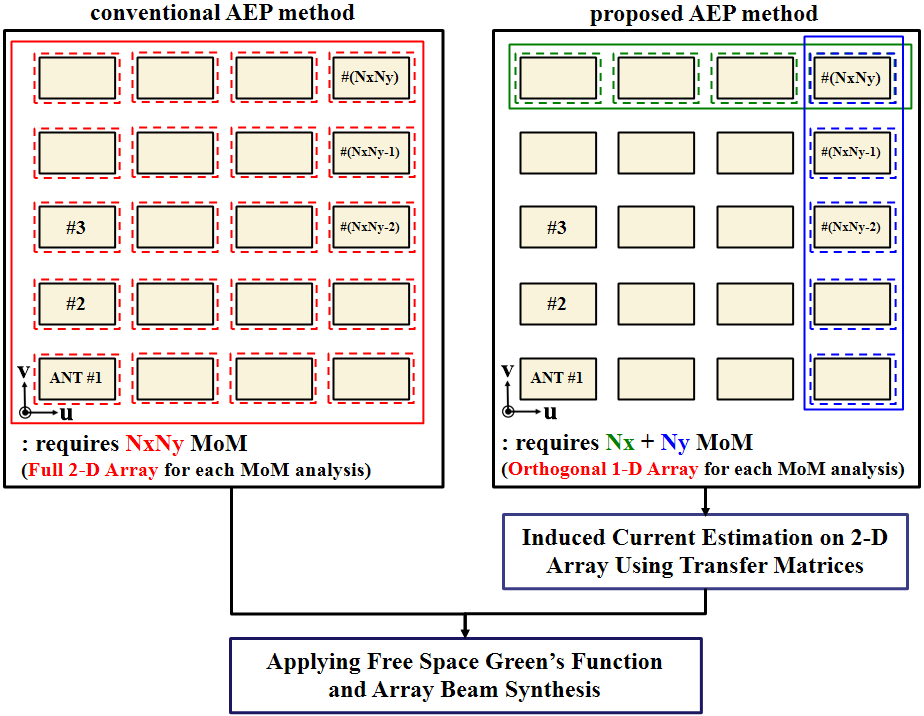}
 \caption{Flowchart of obtaining AEP for the conventional full-wave method and the proposed method}\label{fig1}
\end{figure}
For an isolated single antenna segmented into $M$ triangular meshes, the current density distribution is denoted as $[\mathbf{J}_{\text{iso}}] \in \mathbb{C}^{M \times M}$. In array configurations, port excitation induces currents in terminated elements through electromagnetic coupling. This energy transfer mechanism is represented by the current-based transfer matrix $[\mathbf{C}_k] = [\mathbf{C}_{1,k}, \dots, \mathbf{C}_{i,k}, \dots, \mathbf{C}_{N,k}] \in \mathbb{C}^{M \times N}$, where $\mathbf{C}_{i,k} = [C^{(1)}_{i,k}, \dots, C^{(m)}_{i,k}, \dots, C^{(M)}_{i,k}]^\mathrm{T}$ indicates coupling coefficients between elements.
The transfer matrix $[\mathbf{C}_k$] and current distribution $[\mathbf{J}_k] \in \mathbb{C}^{M \times N}$ on the array antenna for $k$-th port excitation are related as
\begin{equation}
\Big[\mathbf{C}_k\Big] = \Big[\mathbf{J}_{\text{iso}}\Big]^{-1} \Big[\mathbf{J}_k\Big]
\end{equation}
For mesh-to-mesh coupling analysis as shown in Fig. 2 (b), the transfer matrix of the $m$-th mesh is defined as $[\mathbf{C}^{(m)}] = \{ \mathbf{C}^{(m)}_k \}_{k=1}^{N} \in \mathbb{C}^{N \times N}$, characterizing the coupling between excitation ports and array elements.
The current-based transfer matrices are obtained by normalizing 1-D array current distributions with respect to the isolated element case, capturing the coupling effects between excited and terminated ports.
\vspace{-0.2cm}
\subsection{Determination of Transfer Matrix on a 2-D Array}
Current distributions in 2-D array structures can be predicted using transfer matrices derived from orthogonal 1-D array analyses along the $u$-axis and $v$-axis directions. Through the current distributions for each excitation port of the $u$-axis and $v$-axis 1-D array, transfer matrices \([\mathbf{C}_k]|_{\mathbf{u}}\) and \([\mathbf{C}_k]|_{\mathbf{v}}\) is determined using (3) for each axis $u$ and $v$, respectively.
From the determined \([\mathbf{C}^{m}]|_{\mathbf{u}}\) and \([\mathbf{C}^{m}]|_{\mathbf{u}}\), the transfer matrix for the complete 2-D array structure $[\mathbf{C}^m]|_{\mathbf{2D}} \in \mathbb{C}^{N_{x}N_{y} \times N_{x}N_{y}}$ is obtained as
\begin{equation}
\Big[\mathbf{C}^{m}\Big]\Big|_{\mathbf{2D}} = \Big[\mathbf{C}^{m}\Big]\Big|_{\mathbf{u}} \otimes \Big[\mathbf{C}^{m}\Big]\Big|_{\mathbf{v}}
\end{equation}
where $\otimes$ represents the Kronecker product, assuming independent mutual coupling effects in orthogonal directions.
For $k$-th port excitation in the 2-D array, the estimated current distribution is predicted as
\begin{equation}
\Big[\mathbf{J}^{\text{est}}_k\Big] = \Big[\mathbf{J}_{\text{iso}}\Big] \Big[\mathbf{C}_k\Big]\Big|_{\mathbf{2D}}
\end{equation}
where $[\mathbf{C}_k]|_{\mathbf{2D}} = \{ [\mathbf{C}_k^{(m)}]|_{\mathbf{2D}} \}_{m=1}^{M}$.
The radiation pattern at arbitrary observation points is calculated through the free-space dyadic Green's function as
\begin{equation}
E^{\text{AEP}}_{k}(\theta, \phi) = \sum_{m=1}^{M_{total}} G(\mathbf{u}(\theta, \phi), \mathbf{r}_m) \cdot \mathbf{J}^{\text{est}}_{k}(\mathbf{r}_m)
\end{equation}
where $\mathbf{r}_m$ indicates the $m$-th source point vector and $M_{total}(=MN_{x}N_{y})$ represents the total number of mesh elements.

The directional decomposition approach significantly reduces computational complexity while preserving both electromagnetic coupling effect and truncation effects. This efficiency becomes particularly advantageous for highly-integrated antenna structures such as phased arrays and reconfigurable intelligent surfaces. Through current-based formulation and Green's function, the method enables various radiated field synthesis problems [19].

\section{Numerical Validation for Proposed Method}
The effectiveness of the proposed directional decomposition approach is required to be validated in cases of large array structures where computational efficiency becomes critical.
For numerical validation, an open-ended square waveguide antenna with dipole feeding [20-23] is selected due to its strong mutual coupling characteristics in densely arranged rectangular array configurations [20-21].
The array antenna consists of an 11$\times$9 element 2-D array with interelement spacings of 0.14$\lambda$ and 0.12$\lambda$ along $x$ and $y$ axes, respectively, operating at 10GHz, as illustrated in Fig. 2 (a).
\begin{figure} [h]
  \centering
  \subfigure[]{
  \includegraphics[width=0.48\textwidth]{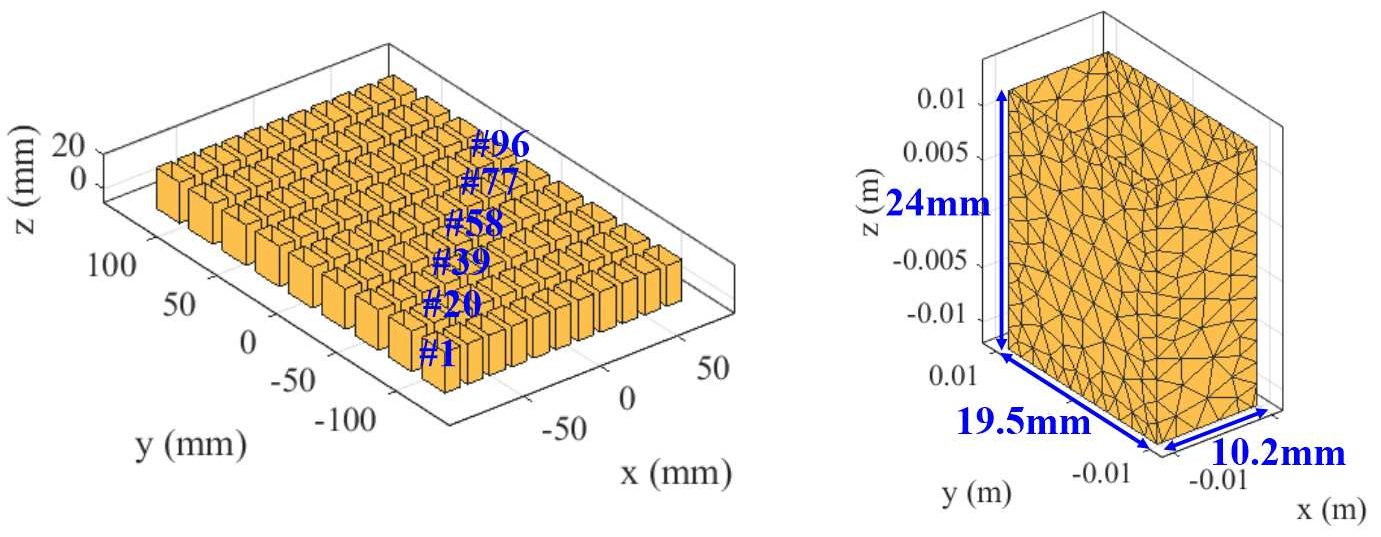}}
  \subfigure[]{
  \includegraphics[width=0.45\textwidth]{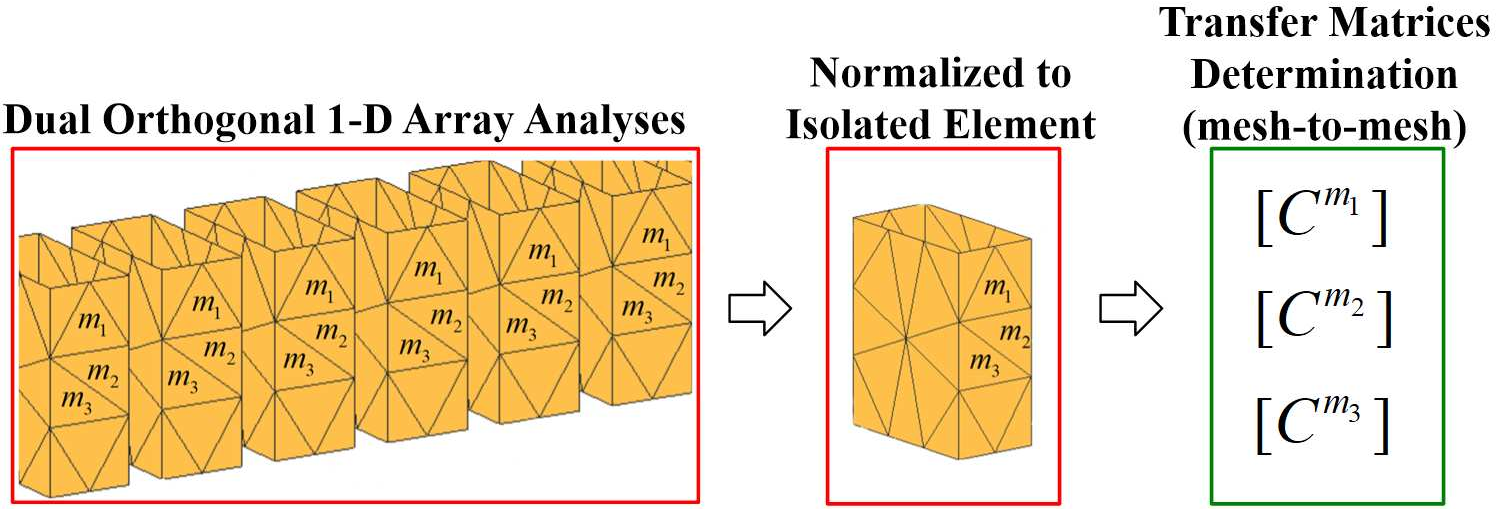}}
  \vspace{-0.2cm}
  \caption{(a) Structure and geometry of 11 $\times$ 9 open-ended square waveguides with port indices (b) Expression of mesh-to-mesh transfer matrices determination}\label{fig2}
\end{figure}

\vspace{-0.2cm}
\subsection{Obtaining Current Distribution on 2-D Array}
The current distribution prediction for 2-D array structures requires full-wave analyses of orthogonal 1-D arrays. The current distributions $[\mathbf{J}_{k}]|_{\mathbf{u}}$ and $[\mathbf{J}_{k}]|_{\mathbf{v}}$ are obtained through MoM simulations supported by FEKO, where each array element is discretized into 598 triangular meshes. For the array along the $x$-axis ($u$ = 1 to 11) and $y$-axis ($v$ = 1 to 9), the analysis requires 6,578 (= 598 × 11) and 5,382 (= 598 × 9) total mesh elements, respectively. Transfer matrices are then constructed using (3) and (4), enabling prediction of current distributions for all 99 antenna ports in the complete 2-D array structure, using (5).
\begin{figure}[h]
 \centering
 \subfigure[]{
   \includegraphics[width=0.45\linewidth]{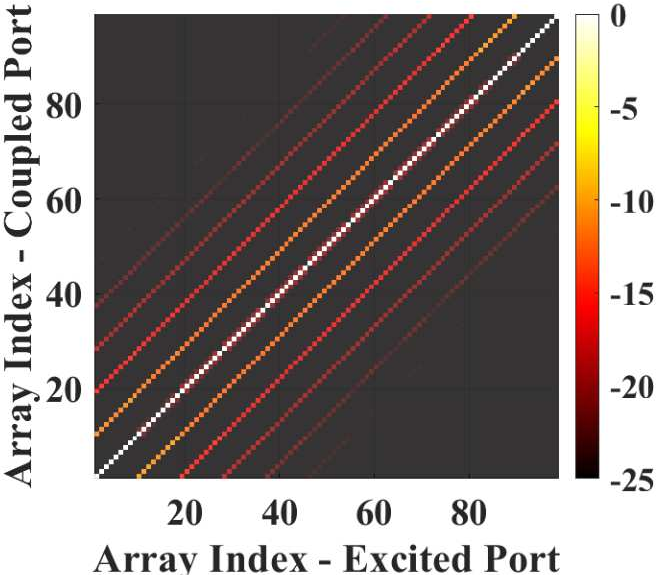}
 }
 \subfigure[]{
   \includegraphics[width=0.45\linewidth]{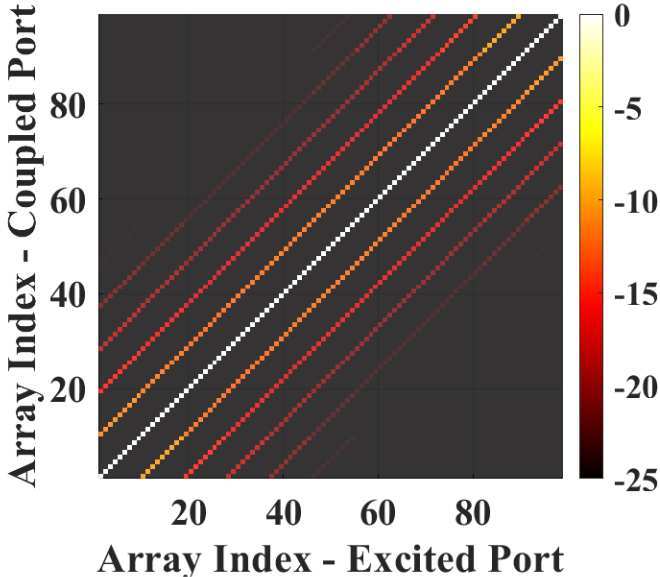}
 }

 \caption{Averaged current distribution spectrum for 99 cases of unit excitation. (a) MoM results (reference). (b) Predicted results (proposed).}
 \label{fig3}
\end{figure}
\begin{figure*} [th!]
 \centering
 \subfigure[]{
   \includegraphics[width=0.99\linewidth]{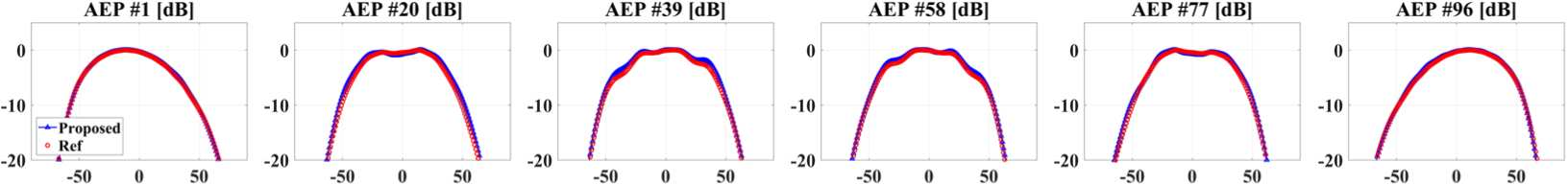}}
 \subfigure[]{
   \includegraphics[width=0.99\linewidth]{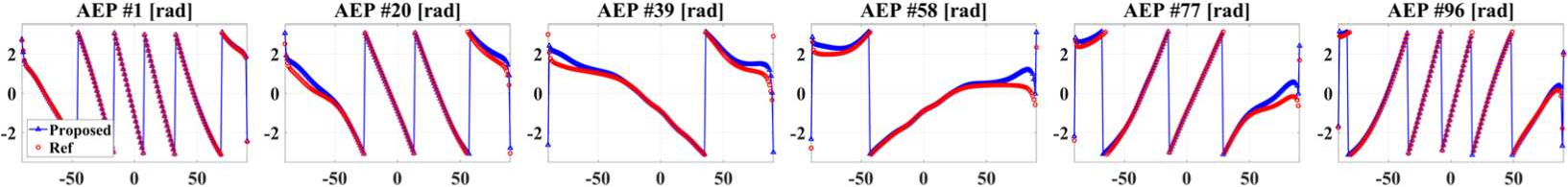}}
   \vspace{-0.2cm}
 \caption{AEP obtained with proposed method (blue line) and reference obtained with MoM (red dot) for $\phi=0^{\circ}$, $\theta=-90^{\circ}$ to $90^{\circ}$ in (a) normalized magnitude pattern and (b) phase pattern.}\label{fig4}
\end{figure*}

In contrast, conventional full-wave AEP method necessitates 99 independent MoM simulations, requiring 59,202 (= 598 × 99) mesh elements for the entire array structure. Fig. 3 presents a current distribution spectrum comparison, where each pixel represents the averaged magnitude of 598 complex-valued currents per array element in logarithmic scale. The current spectrum visualization enables direct observation of mutual coupling effects between antenna elements, particularly demonstrating the strong coupling characteristics along the $x$-axis observed in Fig. 3. The nearly identical distributions between conventional method and predicted results validate the accurate representation of coupling effects through the proposed directional decomposition approach.

\subsection{AEP Estimation}
The radiation pattern calculation utilizes the 99 predicted current distributions with applying Green's function in (8). Fig. 4 presents a comparison of the estimated Apng (proposed) with MoM results (reference) for six uniformly selected antenna ports, indicated in Fig. 2 (a). The E-polarization field pattern, dominant in this array antenna configuration, is compared in the plane of $\phi=0^{\circ}$ for $\theta=-90^{\circ}$ to $90^{\circ}$.
\begin{figure} [t]
 \centering
 \subfigure[]{
 \includegraphics[width=0.99\linewidth]{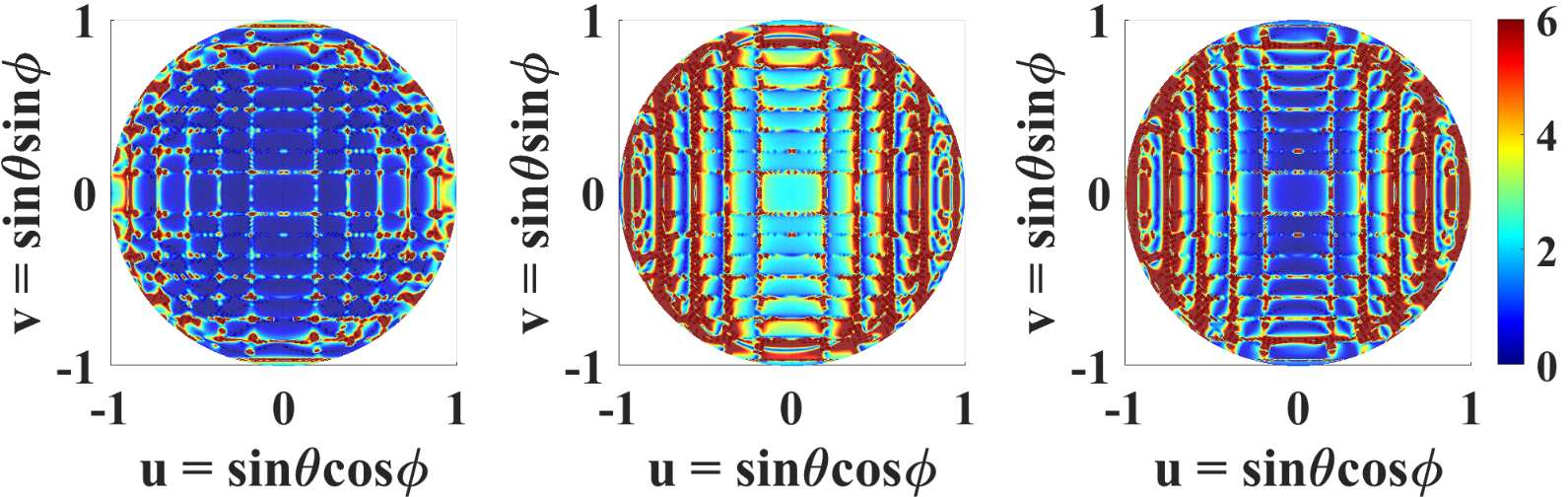}}
 \subfigure[]{
 \includegraphics[width=0.99\linewidth]{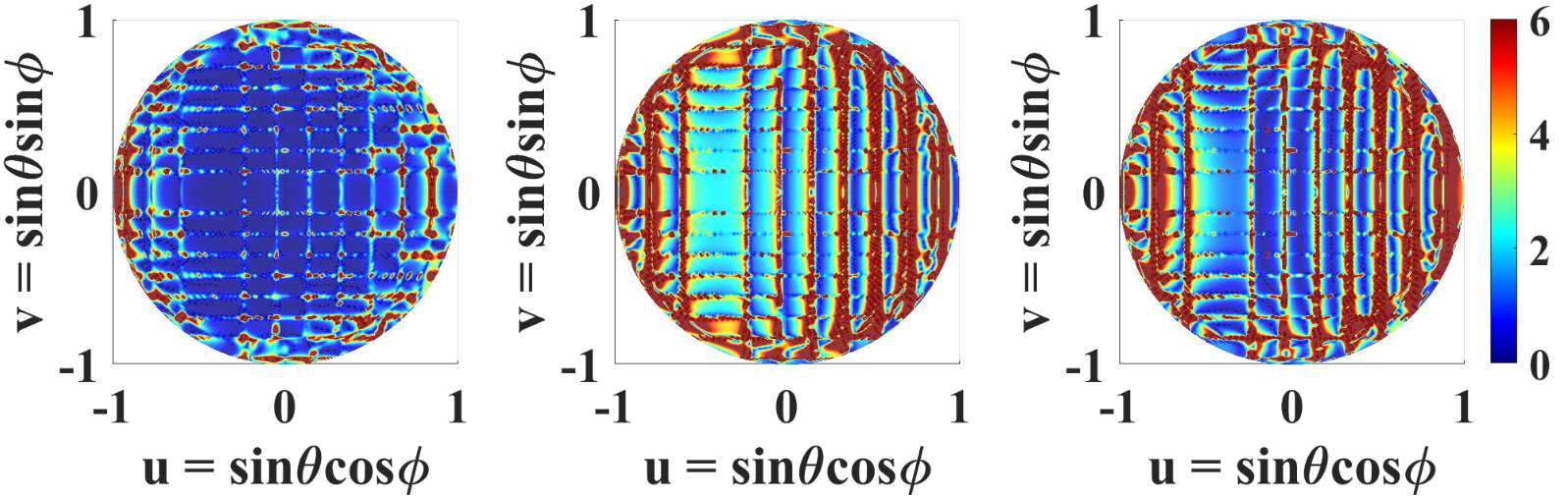}}
 \vspace{-0.2cm}
 \caption{Prediction logarithmic error of synthesized radiation pattern between MoM and PMM-isolated (left), PMM-periodic (center), proposed method (right) for each steering case (a) $\theta=0^{\circ}$, $\phi=0^{\circ}$ and (b) $\theta=25^{\circ}$, $\phi=180^{\circ}$.}\label{fig5}
\end{figure}
The proposed method demonstrates accurate prediction of both magnitude and phase patterns for each antenna port.
Results in Fig. 4 can lead to the reliability of subsequent beam synthesis problems of the array antenna.

\subsection{Beam Synthesis Problem}
For large array structures where conventional AEP method becomes computationally intensive, PMM often serves as an alternative approach [4-8].
While PMM-isolated case utilizing single element patterns neglects both mutual coupling and truncation effects, PMM-periodic case assuming an infinite array condition partially accounts for mutual coupling but fails to capture truncation effects.
In contrast, the proposed method accurately estimates Apng by incorporating both effects, resulting in significantly reduced prediction errors.
Fig. 5 compares the logarithmic radiation pattern errors in $u$-$v$ space between the PMM-based and the proposed beam synthesis results for two steering cases: ($\theta = 0^\circ$, $\phi = 0^\circ$) and ($\theta = 25^\circ$, $\phi = 180^\circ$), as shown in Fig. 5 (a) and Fig. 5 (b), respectively.

\begin{figure} [b]
 \centering
 \includegraphics[width=0.47\textwidth]{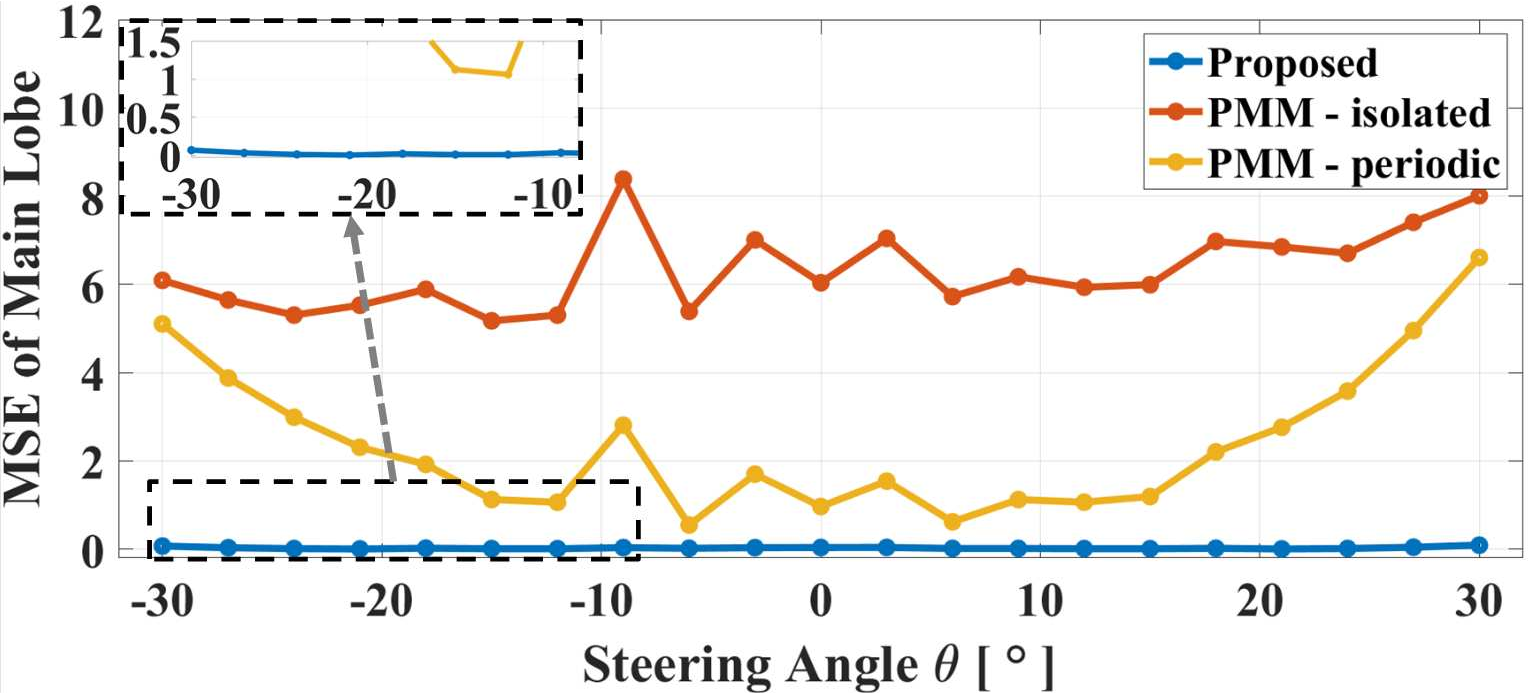}
 \vspace{-0.3cm}
 \caption{Mean square error of main lobe when steering in elevation direction ($\theta=-30^{\circ}$ to $30^{\circ}$) with convex optimization.}\label{fig7}
\end{figure}
Fig. 6 presents synthesized beam patterns through various optimization algorithms, demonstrating the applicability of the proposed method.
The array excitation coefficients were determined through convex optimization [6] to maximize directional gain while minimizing sidelobe levels for steering angles from $\theta = -30^{\circ}$ to $30^{\circ}$ at $\phi = 0^{\circ}$.
Using these optimized coefficients, beam patterns were synthesized with element patterns from PMM-isolated, PMM-periodic, and the proposed method, respectively.
Fig. 6 shows the mean squared error (MSE) in the main lobe region with respect to full-wave analysis results across different steering angles.
The proposed method maintains MSE below 0.1dB for all steering angles, while both PMM approaches show increasing errors at larger steering angles due to intensified mutual coupling and truncation effects.

The computational efficiency between the proposed method and conventional method in AEP extraction is compared in Table 1.
Analyzing two mutually orthogonal 1-D array structures requires significantly lower computational complexity, reducing the calculation time to 0.095\% of that required for full-wave analysis.
This dimensional reduction in computational complexity becomes increasingly advantageous as array structures grow in size and complexity.
\begin{table}[h!]
\caption{Computation Efficiency Comparison Between Proposed and Conventional AEP Methods}
\centering
\begin{tabular}{c|c|c}
\hline\hline  % 상단 이중 실선
& {Proposed} & {Conventional} \\ \hline
\makecell{Number of \\ Triangular Meshes} & $6,578 + 5,382^*$ & $59,202$ \\ \cline{1-3}
\makecell{Computational \\ Complexity$^{**}$} & $\mathcal{O}(M_B^2(N_x^3 + N_y^3))$ & $\mathcal{O}(M_B^2(N_x^3 N_y^3))$ \\ \cline{1-3}
\makecell{Calculation \\ time} & 1.97 min & 34.506 hr \\ \hline\hline  % 하단 이중 실선
\end{tabular}

\vspace{0.15cm}

\hfill \(M_{B}\) : number of basis function for single antenna \\
\hfill * analyzed separately on each orthogonal axis \\
\hfill ** except additional factorization \\
\label{table_comparison}
\end{table}
% ===============================================================================
%% 5. Conclusion
% \vspace{-0.3cm}
\section{Conclusion}
 This letter have presented efficient analysis method for
 obtaining AEP of 2-D array antenna. Using currend-based
 transfer matrix, current distribution of the entire planar array
 for each excitation ports can be predicted, and AEP reflecting
 both mutual coupling and truncation effects can be obtained
 through Green’s function. The AEP estimated by the proposed
 method allows accurate beam synthesis to be performed, which
 is verified through comparison to MoM results. The proposed
 method is applicable not only to open-ended waveguides but
 also to antennas with more complex structures. As future work,
 it can be expected to obtain the AEP of complex array antennas
 with dielectric structures

\vfill

\end{document}